\documentclass[12pt]{amsart}
\usepackage{amsmath,amsthm,amsfonts,amssymb}

\newcount\casecount
\newcommand{\case}[1]{\advance\casecount by 1\vskip .1in
	\goodbreak\noindent\textbf{Case \number\casecount:}\ #1\par\noindent}

\newcommand{\N}{\mathbb{N}}

\newcommand{\SO}{\mathrm{SO}}

\newcommand{\SU}{\mathrm{SU}}

\newcommand{\GL}{\mathrm{GL}}

\newcommand{\SL}{\mathrm{SL}}

\newcommand{\Sp}{\mathrm{Sp}}

\newcommand{\GO}{\mathrm{GO}}

\newtheorem{thm}{Theorem}

\newtheorem{cor}[thm]{Corollary}

\newtheorem{prop}[thm]{Proposition}
\newtheorem{lemma}[thm]{Lemma}

\newtheorem{lem}[thm]{Lemma}

\newtheorem*{claim*}{Claim}
\begin{document}

\title{Low Degree Representations of Simple Lie Groups}

\author{Robert Guralnick}
\address{Department of Mathematics
University of Southern California
3620 Vermont Avenue
Los Angeles, California 90089-2532}
\email{guralnic@usc.edu}

\author{Michael Larsen}
\address{Department of Mathematics,
Indiana University,
Bloomington, IN
U.S.A. 47405}
\email{mjlarsen@indiana.edu}
\thanks{The authors were partially supported by NSF Grants DMS-0653873 and  DMS-0800705.}

\author{Corey Manack}
\address{Department of Mathematics,
Indiana University,
Bloomington, IN
U.S.A. 47405}
\email{cmanack@indiana.edu}
\maketitle

\section{Introduction}

We say a compact Lie group is \emph{simple} if it is connected,  has finite center
and is a simple group modulo its center.  
This paper is mostly concerned with upper bounds for $R_n(G)$, 
the number of characters of $G$ of degree 
$\le n$.  If is convenient to define the \emph{Witten zeta-function}
$$\zeta_G(s) = \sum_{\chi}\chi(1)^{-s},$$
where the sum is taken over all irreducible characters of $G$.  By \cite{LL}, the abscissa of convergence of this 
Dirichlet series is $2/h$, where $h$ denotes the Coxeter number of $G$.
This implies
$$R_n(G) = O(n^{2/h+\epsilon})$$
for all $\epsilon > 0$.   One can turn this around and try to bound $R_n(G)$ for fixed $n$ and varying $G$.  In this direction, Wallach \cite{nw} proved that $R_n(G) \le 7n$.

Our first result gives the optimal bound for this problem.

\begin{thm}
\label{guralnick}
If $G$ is a simple compact Lie group and $n$ is a positive integer,
$$R_n(G)\le n$$
with equality only if $n=1$, $G=\SU(2)$, or $G=\SU(3)$ and $n=3$.
\end{thm}

If we exclude groups of low Coxeter number, we can improve on Theorem~\ref{guralnick}.
A formulation in terms of the Witten zeta-function is as follows.

\begin{thm}
\label{limit}
If $G_1,G_2,\ldots$ is a sequence of pairwise non-isomorphic simple compact Lie groups, and $s>0$, then
$$\lim_{n\to\infty} \zeta_{G_n}(s) = 1.$$
\end{thm}

This can be regarded as an analogue for compact Lie groups of the main theorem of 
\cite{LS}.  An immediate consequence of Theorem~\ref{limit} is

\begin{cor}
\label{strong-guralnick}
For every $\epsilon > 0$ there exists $N$ such that if $G$ is a simple compact Lie group of
dimension $\ge N$, then $R_n(G) \le n^\epsilon$ for all $n\ge 1$.
\end{cor}

As a further application, we prove the following theorem.

\begin{thm}
\label{repeats}
For all $\epsilon > 0$, the number of isomorphism classes of pairs $(G,V)$, where $G$ is a 
compact semisimple group and $V$ is a faithful, irreducible representation of dimension $n$ is 
$O(n^\epsilon)$.
\end{thm}

If $G$ is either a connected algebraic group or a compact Lie group, let $\mathcal{M}(G)$
denote the set of conjugacy classes of maximal closed subgroups of $G$ (i.e. maximal
subgroups among closed proper subgroups). 
The previous result can be used to give bounds for $|\mathcal{M}(G)|$ in
characteristic $0$.   Considering parabolic
subgroups (or subspace stabilizers for classical compact groups), we see that $|\mathcal{M}(G)|$
is greater than the rank of $G$.   Using the previous result, we can prove:

\begin{thm} 
\label{maximal}  Let $G$ be a simple algebraic group over an algebraically
closed field of characteristic $0$  or a compact Lie group of rank $r$.
Then  $|\mathcal{M}(G)|=O(r)$.
\end{thm}

With some effort, one could get a more precise upper bound in the previous result. 
The previous result uses a result of H\"as\"a \cite{hasa} which does depend upon
the classification of finite simple groups.  The result for positive dimensional maximal
closed subgroups does not depend on this. 

The last result of the paper is a result about dimensions of fixed spaces of elements in
compact groups.    The first author and Malle \cite{gm}  
showed  that if $G$ is an irreducible subgroup of $\GL(n,k), n > 1$ (for any field $k$), then
there is an element $g \in G$ with fixed space of dimension at most $n/3$  (the example
$\SO(3)$ shows that this is best possible).   This answered a conjecture of Peter Neumann.
We show that as $n$ increases we can do much better for Lie groups.   
Let $w(g)=w_V(g)$ denote the largest dimension of any eigenspace of $g$ acting on $V$.

\begin{thm}
\label{gur-malle}   Let $G$ be a simple compact Lie group.   Let $V$ be a finite dimensional
irreducible module.   Let $\epsilon > 0$.   There exists $N = N(\epsilon)$ such that if
$\dim V > N$,  then for a  dense open subset  $X$ of $G$,  $g \in X$ implies
that $w_V(g)  < \epsilon \dim V$.
\end{thm}

This immediately implies the same result holds for algebraic groups over algebraically 
closed fields of characteristic $0$.  

\section{Counting Representations}

We begin by noting that if $G$ is a simple compact Lie group, every representation of $G$ can be
regarded as a representation of its universal covering group $\tilde{G}$.  Thus,
$$R_n(\tilde{G}) \ge R_n(G)$$
for all $n$, and
$$\zeta_{\tilde{G}}(s) \ge \zeta_G(s)$$
for all $s>0$ for which the left hand side is defined.   Consequently, it suffices to prove the two theorems for 
simply connected simple compact Lie groups, and we will assume henceforth that
$G = \tilde{G}$.  This is immediate for Theorem~\ref{limit}; for Theorem~\ref{guralnick}, we note
that if $G$ is a non-trivial quotient of $\SU(2)$ (resp. $\SU(3)$), $R_n(G) < n$ for $n\ge 2$ (resp. for $n=3$).  
We therefore consider only simply connected groups henceforth.

The irreducible characters $\chi$ of $G$ are indexed by dominant weights of $G$.  Let
$r$ be the rank of $G$, and $\varpi_1,\ldots,\varpi_r$ the fundamental weights.  Then the
dominant weights are the non-negative linear combinations 
$$\lambda = (c_1-1)\varpi_1+\cdots+(c_r-1)\varpi_r,$$
where the $c_i$ are positive integers.  The half sum $\rho$ of positive roots equals the sum of the 
fundamental weights \cite[VI, \S1, Prop. 29]{B}, so 
$$c_1\varpi_1 + \cdots + c_r\varpi_r = \lambda+\rho.$$
Weyl's dimension formula for the character $\chi_\lambda$ of highest weight $\lambda$ asserts
$$\chi_\lambda(1) = \prod_{\alpha\in \Phi^+} \frac{\alpha^\vee(\lambda+\rho)}{\alpha^\vee(\rho)},$$
where $\Phi^+$ denotes the set of positive roots of $G$, and $\alpha^\vee$ denotes the dual root of $\alpha$.  
We can express this in terms of the $c_i$ as follows.  If $\alpha_1^\vee,\ldots\alpha_r^\vee$ are the simple positive 
roots of the dual root system $\Phi^\vee$, and
$\phi_1^\vee,\ldots,\phi_u^\vee$ are all the positive roots in $\Phi^\vee$, then every 
root $\beta_i^\vee\in \Phi^\vee$
can be written uniquely as
$$\beta_i^\vee = \sum_{j=1}^r a_{i,j}\alpha_j^\vee,$$
where the $a_{i,j}$ are non-negative integers.  
As $\alpha_i^\vee(\varpi_j) = \delta_{ij}$ \cite[VI, \S1.10]{B}, we have
$$\chi_\lambda(1) = \prod_{i=1}^u \frac{\sum_{j=1}^r a_{i,j} c_j}{\sum_{j=1}^r a_{i,j}}
= \prod_{i=1}^u \sum_{j=1}^r w_{i,j} c_j,$$
where the 
$$w_{i,j} = \frac{a_{i,j}}{\sum_{j=1}^r a_{i,j}}$$
are non-negative and sum to $1$.

By the weighted arithmetic-geometric mean inequality,
$$\sum_{j=1}^r w_{i,j} c_j \ge \prod_{j=1}^r c_j^{w_{i,j}}.$$
Multiplying over $i=1,2,\ldots,u$, we obtain
\begin{equation}
\label{dim-est}
\chi_\lambda(1) \ge \prod_{i=1}^u c_j^{v_j},
\end{equation}
where
\begin{equation}
\label{v-def}
v_j = v_j(\Phi^\vee) =  \sum_{i=1}^u w_{i,j}.
\end{equation}

Defining
$$Z_{\Phi}(s) = \prod_{i=1}^r \zeta(v_j(\Phi)s),$$
where $\zeta(s)$ denotes the Riemann zeta-function, we obtain
\begin{equation}
\label{main-est}
\zeta_G(s) \le Z_{\Phi^\vee}(s)
\end{equation}
for all $s$ for which the right hand side converges.

Our next task is to estimate values of $Z_{\Phi}(s)$ for the various root systems
$\Phi$.

\begin{lem}
\label{ABCFG}
If $\Phi$ is of type $A_r$ and the simple roots $\alpha_i$ are numbered as usual, then
$$v_j = \sum_{i=1}^j\sum_{k=j}^r \frac 1{1+k-i}.$$
\end{lem}

\begin{proof}
The roots of $\Phi$ are precisely the sums of consecutive simple roots,
$$\alpha_i + \alpha_{i+1} + \cdots + \alpha_k.$$
Each such root contributes $\frac 1{1+k-i}$ to the sum (\ref{v-def}) if $i\le j\le k$ and $0$ otherwise.
\end{proof}

\begin{lem}
\label{v-est}
If $\Phi$ is of type A and $j\le \frac{r+1}2$, then
\begin{equation}
\label{ineq-1}
v_j = v_{r+1-j} \ge j\max(1,\log r-\log j).
\end{equation}
\end{lem}

\begin{proof}
The equality $v_j = v_{r+1-j}$ is clear from Lemma~\ref{ABCFG}.
If $1\le m\le n\le j$, then $1\le 1+n-m \le j \le j+n-m\le r$.
Setting $i=1+n-m$ and $k = j+n-m$, we obtain
$$v_j \ge \sum_{i=1}^j\sum_{k=j}^r \frac 1{1+k-i} \ge \sum_{n=1}^j\sum_{m=1}^n \frac 1n = j.$$
On the other hand,
$$v_j \ge \sum_{i=1}^j\sum_{k=j}^r \frac 1{1+k-i} \ge \sum_{i=1}^j \sum_{k=j}^r \frac 1{1+k-i}
\ge j\sum_{k=j}^r \frac 1k \ge  j\int_j^r \frac {dx}x.$$
\end{proof}

\begin{lem}
\label{zeta-est}
For $s\ge 3$,
$$\log \zeta(s) \le 2^{1-s}.$$
\end{lem}

\begin{proof}
As $x^{-s}$ is convex for $x>0$,
$$\zeta(s) < 1^{-s}+2^{-s} + \int_{5/2}^\infty x^{-s}\,dx = 1+2^{-s} + \frac{(5/2)^{1-s}}{s-1}.$$
As $s\ge 3$,
$$(5/2)(5/2)^{-s} \le (5/2) (5/4)^{-3} 2^{-s} < \frac{32}{25}\cdot 2^{-s},$$
so 
$$\zeta(s) < 1+2^{-s} + \frac{2^{1-s}}{2} < 1+ 2\cdot 2^{-s},$$
which implies the lemma.
\end{proof}

\begin{prop}
\label{eval-at-1}
If $\Phi$ is of type $A_r$ with $r\ge 9$, of type $D_r$ with $r\ge 5$, or of type $E_r$, then
$$\zeta_G(1) < 13/7.$$
\end{prop}

\begin{proof}
Table 1 gives a direct machine computation of the values of $Z_\Phi(1) = Z_{\Phi^\vee}(1)$
to four decimal places.

We begin with type $A_r$.
Table 1 shows that this inequality is true for $9\le r\le 20$.
We therefore assume $r\ge 21$.
By (\ref{main-est}), Lemma~\ref{v-est}, and Lemma~\ref{zeta-est},
\begin{align*}
\zeta_G(1) &\le \prod_{i=1}^6 \zeta(v_i)^2 \prod_{i=7}^{(r+1)/2} \zeta(v_i)^2 
\le \prod_{i=1}^6 \zeta(i \log(21/i))^2  \prod_{i=7}^{(r+1)/2} \zeta(v_i)^2 \\
& \le \prod_{i=1}^6 \zeta(i \log(21/i))^2  \prod_{i=7}^\infty  e^{2^{2-i}} 
 = \prod_{i=1}^6 \zeta(i \log(21/i))^2 e^{1/16} < 13/7.
\end{align*}

Next, we consider type $D_r$, using the standard numbering of simple roots \cite[Planche IV]{B}.
If $1\le i\le r-2$ and $1\le j\le r$, there is a unique path from the $i$th to the $j$th vertex in 
the Dynkin diagram; its length is $|j-i|$ unless $j=r$, in which case it is $j-1-i = |j-i|-1$.
By \cite[VI, \S1, Prop. 19, Cor. 3]{B}, the sum of the simple roots corresponding to the vertices in
this path is again a positive root.  It follows that for $1\le i\le r$, 
$$v_i(D_r) \ge v_i(A_r).$$
If $1\le p < q \le r-1$, then
\begin{equation}
\label{D-roots}
\sum_{i=p}^{r} \alpha_i + \sum_{j=q}^{r-2} \alpha_j
\end{equation}
is a positive root.  So is
$$\beta + \sum_{i=p}^{r-2} \alpha_i,$$
for $\beta\in \{\alpha_{r-1},\alpha_r\}$ and $1\le p\le r-1$.  Thus,
$$v_{r-1}(D_r) = v_r(D_r) \ge \sum_{i=1}^{r-1} \frac{\lceil (i+1)/2\rceil}{i} \ge \frac{r-1}2.$$
If $r\ge 11$, then 
$$v_{r-1}(D_r) > v_{r-1}(A_r) > v_r(A_r).$$
Therefore, in this range the proposition is implied by the result for type $A_r$.
For $5\le r\le 10$, it follows from Table 1.

Likewise, for $E_6$, $E_7$, and $E_8$, it follows from Table 1.

\begin{table}
\centering
\begin{tabular}{|r|l|r|l|}
\hline
$A_9$&1.8558 & $D_5$&1.7269 \\
$A_{10}$&1.7336 & $D_6$&1.4609 \\
$A_{11}$&1.6400 & $D_7$&1.3244 \\
$A_{12}$&1.5667 & $D_8$&1.2462 \\
$A_{13}$&1.5083 & $D_9$&1.1978 \\
$A_{14}$&1.4610 & $D_{10}$&1.1160 \\
$A_{15}$&1.4221 & $E_6$&1.2522 \\
$A_{16}$&1.3896 & $E_7$&1.0836 \\
$A_{17}$&1.3621 & $E_8$&1.0178 \\
$A_{18}$&1.3386 && \\
$A_{19}$&1.3182 && \\
$A_{20}$&1.3004 && \\
\hline 
\end{tabular}
\medskip
\caption{$Z_\Phi(1)$}
\end{table}

\end{proof}

We can now prove Theorem~\ref{guralnick}.

\begin{proof}
If $\Phi$ is among the root systems covered by Proposition~\ref{eval-at-1}, then the same is true for its dual, 
and it follows that $\zeta_G(1) < 13/7$.  If, for some $n\ge 7$,
we have $R_n(G) \ge n$, then
$$\zeta_G(1) > 1 + (n-1)\cdot n^{-1}\ge 13/7.$$
Thus, $n\le 6$.  However, $\dim G \ge 45 > 6^2$ for all the cases 
covered by Proposition~\ref{eval-at-1},
so $G$ has no non-trivial representations of degree $\le 6$.

Table 2 indicates the values of $\zeta_G(3/4)^4$ for all simply laced root systems of
rank $\ge 2$ not covered by Proposition~\ref{eval-at-1}.
An exhaustive machine search for representations in the indicated ranges yields precisely one
non-trivial case where $R_n(G) \ge n$, namely the equality $R_3(\SU(3)) = 3.$
Theorem~\ref{guralnick} follows for all $G$ with simply laced root system.

The remaining irreducible root systems all have the same underlying graph as $A_r$, and we use the same 
numbering convention for the simple roots.
By \cite[VI, \S1, Prop. 19, Cor. 3]{B}, $\alpha_i+\alpha_{i+1}+\cdots+\alpha_j$ is a positive root
for all $i\le j$.  Therefore, the dimension of the irreducible representation of $G$ associated to
the $r$-tuple $(c_1,\ldots,c_r)$ is at least as great as the dimension of the irreducible representation of $\SU(r+1)$ 
associated to the same $r$-tuple.  
This proves $R_n(G)\le n$ for root systems of type B, C, F, and G as well.
Moreover, in the case $r=2$,
the pairs $(1,2)$ and $(2,1)$ give representations of $B_2$ and $G_2$ of dimension strictly greater than $3$.  
Thus, the non-simply laced root systems give no additional cases for which
$R_n(G) = n>1$. 

\begin{table}
\centering
\begin{tabular}{|r|l|}
\hline
$A_2$&29547000 \\
$A_3$&194100 \\
$A_4$&27475 \\
$A_5$&7055.8 \\
$A_6$&2414.5 \\
$A_7$&1001.9 \\
$A_8$&481.56 \\
$D_4$&994.94\\
\hline 
\end{tabular}
\medskip
\caption{$Z_\Phi(3/4)^4$}
\end{table}

\end{proof}

Now we prove Theorem~\ref{limit}.

\begin{proof}
Since we are taking the limit as $n\to \infty$, we can exclude the exceptional groups and groups
of type $D_r$, $r\le 8$.  Then, if $G$ is of rank $r$ and $s > 2/(r+1)$,
$$\zeta_G(s) \le \zeta_{\SU(r+1)}(s).$$
Thus, it suffices to prove that for all $s>0$.
\begin{equation}
\label{SU-case}
\lim_{r\to\infty} \zeta_{\SU(r+1)}(s) = 1.
\end{equation}
By (\ref{ineq-1}),
$$\zeta_{\SU(r+1)}(s) 
< \prod_{i=1}^{(r+1)/2} \zeta\Bigl(s\sum_{i=1}^j\sum_{k=j}^r \frac 1{1+k-i}\Bigr)^2
< \prod_{j=1}^\infty \zeta(s\max(j,j\log r/j))^2.$$
If we assume $r > (3/s)e^{3/s}$, for each positive integer $j$, either $j > 3/s$ or $r/j > e^{3/s}$, so by
Lemma~\ref{zeta-est},
$$\zeta_{\SU(r+1)}(s) < \exp\Bigl(4\sum_{j=1}^\infty 2^{-s\max(j,j\log r/j)}\Bigr).$$
If $N > \exp(1/s)$ and $r>N^2$,
$$\sum_{j=1}^\infty 2^{-s\max(j,j\log r/j)} \le \sum_{j=1}^N 2^{-sj\log N}
+\sum_{j=N+1}^\infty 2^{-sj} \le 2^{1-s\log N} + \frac{2^{-sN}}{1-2^{-s}}.$$
If $s>0$ is fixed and $N\to\infty$, the right hand side goes to $0$.  As $\zeta_G(s) > 1$ for all $G$
and all $s$, this implies (\ref{SU-case}).

\end{proof}

\begin{lemma}
\label{semi}
Let $a_1,a_2,\ldots$ be a sequence of non-negative real numbers.
If $a_n = O(n^\epsilon)$ for all $n>0$, then
$$b_n := \sum_{\substack{d_1\cdots d_k = n \\ 1 < d_1\le\cdots\le d_k}} = O(n^\epsilon)$$
for all $\epsilon > 0$.
\end{lemma}

\begin{proof}
If
$$f(s) := 1^{-s} + \sum_{n=2}^\infty a_n n^{-s}$$
converges for all $s>0$, then 
$$\prod_{n=2}^\infty \frac{1}{1-a_n n^{-s}} = 1 + \sum_{n=2}^\infty b_n n^{-s}$$
converges for all $s>0$.   If $a_n = O(n^\epsilon)$ for all $\epsilon > 0$, then $f(s)$
converges for all $s>0$, so $b_n = O(n^\epsilon)$ for all $\epsilon > 0 $.
\end{proof}

\begin{lemma}
\label{bounded-weights}
For all positive integers $n$ there exist integers $R$ and $m$ and a finite set 
$\Lambda \subset \N^m$
such that if $G$ is a simply connected, simple, compact Lie group of rank $r\ge R$
and $\chi_\lambda(1) \le r^n$, then either $G$ is of type $A_r$ and
\begin{equation}
\label{A-only}
\lambda = a_1\varpi_1+\cdots+a_m\varpi_m + b_m\varpi_{r+1-m} + \cdots + b_1 \varpi_r
\end{equation}
for some $(a_1,\ldots,a_m),(b_1,\ldots,b_m)\in \Lambda$, or
$G$ is of type $B_r$, $C_r$, or $D_r$, and
\begin{equation}
\label{BCD}
\lambda = a_1\varpi_1+\cdots+a_m\varpi_m
\end{equation}
for some $(a_1,\ldots,a_m)\in \Lambda$.
\end{lemma}

\begin{proof}
For $1\le i\le r$, 
$$v_i(\Phi)\ge 1+1/2+\cdots+1/{r-1} > \log r.$$
By (\ref{dim-est}), this implies $c_i \le m/\log 2$ for all $i$.  Our goal is therefore, to prove
that $c_i = 1$ except when $i$ is bounded above or, in the case of $A_r$, when $r-i$ is bounded above.  
For $A_r$, this is an easy consequence of Lemma~\ref{ABCFG}.  For the remaining root systems, we need 
only prove that for every constant $C$, $v_i(\Phi) \le C\log r$ implies an upper bound on $i$.  
For $B_r$, by \cite[Planche II]{B}, for $1\le p < q \le r$,
$$\sum_{i=p}^{q-1}\alpha_i + \sum_{j=q}^r \alpha_i \in \Phi^+.$$
Thus, the sum for $v_i(B_r)$ contains at least $i(2r-1-i)/2$ terms, each of which is at least $1/(2r-1)$.  Thus,
$$v_i(B_r) \ge \max(v_i(A_r),\frac{i(r-1)}{2r-1}),$$
so $v_i(B_r) \le C\log r$ implies an upper bound on $i$ independent of $r$.
For $C_r$, by \cite[Planche III]{B}, for $1\le p < q \le r$,
$$\alpha_r+\sum_{i=p}^{q-1}\alpha_i + \sum_{j=q}^{r-1} \alpha_i \in \Phi^+.$$
As in the $B_r$ case, $v_i(C_r)$ is the sum of at least  $i(2r-1-i)/2$ terms, each of which is at least $1/(2r-1)$, so
$$v_i(C_r) \ge \max(v_i(A_r),\frac{i(r-1)}{2r-1}),$$
which implies an upper bound on $i$, independent of $r$.
By (\ref{D-roots}), 
$$v_i(D_r) \ge max(v_i(A_r),frac{i(r-1)}{2r-3}),$$
for $1\le i\le r-2$, and
$$v_{r-1}(D_r) = v_r(D_r) \ge \frac{(r-1)(r-2)}{2r-3},$$
which again implies an upper bound on $i$, independent of $r$.
\end{proof}

We now prove Theorem~\ref{repeats}.

\begin{proof}
If $G$ is a semisimple compact Lie group with universal covering group $\tilde G$ and $V$ is a faithful 
irreducible representation of $G$, then we can regard $V$ as an almost faithful irreducible representation of $\tilde G$.  
Writing $\tilde G$ as a product of its simple factors $G_1\times\cdots\times G_k$, we can factor $V$ as 
$V_1\boxtimes \cdots\boxtimes V_k$, where each 
$V_i$ is a non-trivial representation.  We can order the factors $G_i$ so that
$\dim V_1\le \dim V_2\le\cdots\le \dim V_k$.  By Lemma~\ref{semi}, it suffices to prove 
that for each $n$, the number of pairs $(G,V)$ where $G$ is a simply connected, simple, compact Lie group and 
$V$ is a representation of dimension $n$ is $O(n^\epsilon)$ for all $\epsilon > 0$.

By Theorem~\ref{limit}, for every $\epsilon > 0$, there exists 
a finite set $\Sigma$ of groups such that
for every simply connected, simple, compact Lie group $G$ not isomorphic to any element of
$\Sigma$, we have $\zeta_G(\epsilon/2) \le 2$.  For any such $G$, $R_n(G) \le n^{\epsilon/2}$.
If we limit ourselves to groups $G$ of rank $\le n^{\epsilon/2}$ which do not lie in $\Sigma$,
the total contribution is bounded by $n^{\epsilon}$.

Next we show that for each $G\in \Sigma$, there are $O(n^{\epsilon})$ representations of $G$
of degree $n$.  Suppose that $G$ has rank $r$.  The representations of $G$ are indexed by
$(c_1,\ldots,c_r)\in \N^r$.  It suffices to show that for each $j\in\{1,\ldots,r\}$, 
the number of representations with $c_j > 0$ and degree $n$ is $O(n^{\epsilon})$.
Applying \cite[VI, \S1, Prop. 19, Cor. 3]{B}, for each $k\in \{1,\ldots,r\}$, the sum of all the simple roots 
corresponding to vertices on the path from $j$ to $k$ is a positive root $\beta_k$.
Clearly, $\beta_1,\ldots,\beta_r$ span the root system.  Ordering the positive roots
of $G$ so that $\beta_1,\ldots,\beta_r$ come first, we deduce from Weyl's dimension formula that
$$\prod_{i=1}^r \sum_{j=1}^r a_{i,j} c_j\Bigm| C_G n,$$
where
$$C_G =  \prod_{i=1}^u \sum_{j=1}^r a_{i,j}.$$
Moreover, the linear independence of $\beta_1,\ldots,\beta_r$ implies that the $r$-tuple
$$(d_1,\ldots,d_r) = \Bigl(\sum_{j=1}^r a_{1,j} c_j,\ldots,\sum_{j=1}^r a_{r,j} c_j\Bigr)$$
determines $(c_1,\ldots,c_r)$.  We define
$$d_{r+1} = \frac{C_G n}{d_1\cdots d_n}$$
and then let $e_1,\ldots,e_s$, $s\le r+1$, denote the sequence obtained by ordering
all elements in the set $\{d_1,\ldots,d_{r+1}\}$ which are greater than $1$.  Thus
$$\sum_{d_1\cdots d_r | C_G n} 1 \le (r+1)! 
\sum_{\substack{e_1\cdots e_s = C_G n\\ 1 < e_1\le e_2\le \cdots \le e_s}} 1.$$
By Lemma~\ref{semi}, this sum is $O((C_G n)^\epsilon)$ and therefore $O(n^{\epsilon})$.

This leaves only groups of rank $\ge n^{\epsilon/2}$, and taking $n$ sufficiently large, we may assume that 
only classical groups need be considered.  Setting $m= 2/\epsilon$
in Lemma~\ref{bounded-weights}, there exists a finite set $\Lambda$, so that if $G$ is of type 
$A_r$, $\lambda$ is given by (\ref{A-only}), and if $G$ is of type $B_r$, $C_r$, or $D_r$, it is given by (\ref{BCD}).  
We claim that for each element $(a_1,\ldots,a_m)\in\Lambda$ and each family A, B, C, D, the dimension of the 
representation is a strictly monotone function of $r$.  Perhaps the easiest way is note that the branching rules for 
classical groups \cite[9.14,9.16,9.18]{K} 
imply that for $X\in \{A,B,C,D\}$,
the restriction of the representation of $X_r$ associated to $(a_1,\ldots,a_m)$ to $X_{r-1}$ contains the representation of 
$X_{r-1}$ associated to $(a_1,\ldots,a_m)$ as a proper subrepresentation.
We conclude that the groups of rank $\ge n^{\epsilon/2}$ contribute at most a constant term
$|\Lambda|^2 + 3|\Lambda|$ to the total number of pairs $(G,V)$.
\end{proof}

\section{Maximal Subgroups of Classical Groups}

Let $k$ be an algebraically closed field of characteristic $p \ge 0$.   
Let $G$ be a simple classical group over $k$
(i.e. $G=\SL(n), n > 1, \Sp(2n), n > 1, \SO(2n), n > 3, \SO(2n+1), n > 2, p \ne 2$).   
We now
give Aschbacher's description of proper subgroups in the case of algebraic groups.  
The proof is much easier
in the algebraically closed case.   See  \cite{As, KL, gur, gt,  ls2}.   Let $V$ denote
the natural module for $G$.

The main idea of the proof is as follows.  Let $H$ be a closed proper subgroup
of $G$.  We may assume that $H$ acts irreducibly and primitively on $V$.
In particular, any normal subgroup $N$ of $H$ acts homogeneously on $V$
(i.e. $V$ is a direct sum of isomorphic irreducible modules of $G$).   Choose
$N$ to be a minimal closed normal noncentral subgroup of $H$.  Then $N$
cannot be abelian (since it acts homogeneously and is noncentral).  So
either $N$ is of symplectic type (i.e. essentially a finite extraspecial group) or it is a central
product of quasimple groups.   If $N$ does not act irreducibly, it follows
by Clifford theory, that $G$ preserves a tensor structure on $V$ --- i.e.,
$V \cong V_1 \otimes V_2$ with each $V_i$ a module for the cover of $N$.
In this case, $N$ will be contained in $X \otimes  Y$ where $X$ is a classical group
acting on $V_i$  (e.g., for $G=\SL$,  $N$ is contained in $GL(V_1) \otimes GL(V_2)$,
or if $G=\SO$, $N$ is contained in $\GO(V_1) \times \GO(F_2)$ or 
$\Sp(V_1) \otimes \Sp(V_2)$).     So we may assume that $N$ is irreducible.
If $N$ is of symplectic type, then $H$ is contained in the normalizer of $N$.
If $N$ is a central product of quasisimple groups, then either $N$ is finite or
semisimple.  If there is one more than one term in the central product, 
it is easy to see that $V$
is tensor induced (i.e. $V=V_1 \otimes \ldots \otimes V_m$  with $\dim V_i=d$
for all $i$ and  $H$ embeds in $\GL(d) \wr S_m$-- if $G = SO(d^m)$, then
in fact $H$ embeds $\GO(d) \wr S_m$ and or possibly $\Sp(d) \wr S_m$;
there are similar possibilities for $G =Sp(d^m)$).
 
Thus, Aschbacher's theorem for algebraic groups can be stated as:   

Let $H$ be a proper closed subgroup of  $G$, then  $H$ belongs
to one of the classes $\mathcal{C}_i(G)$ or $\mathcal{S}_i(G)$.  

\bigskip

\begin{table}
\centering
\begin{tabular}{|r|l}
 \hline
$\mathcal{C}_1$  &   $H$ preserves a totally singular or nondegenerate subspace of $V$   \\  \hline
$\mathcal{C}_2$ &  $H$ preserves an additive decomposition of $V$ \\   \hline
$\mathcal{C}_4$ &  $V$ is tensor decomposable for $H$  \\  \hline
$\mathcal{C}_6$  & $H$ normalizes a subgroup of symplectic type  \\   \hline
$\mathcal{C}_7$  & $V$ is tensor induced for $H$ \\   \hline
$\mathcal{C}_8$ &   $V$ is the natural module for a classical subgroup $H$; \\   \hline
$\mathcal{S}_1$ & $H$ is the normalizer of an irreducible  quasisimple finite group $S$ \\ \hline
$\mathcal{S}_2$ & $H$ is the normalizer of an irreducible simple algebraic subgroup $S$  \\  \hline
\end{tabular}
\medskip
\caption{Aschbacher Classes}
\end{table}

\bigskip

Note that we have kept the same numbering as in Aschbacher's theorem for finite
groups but certain cases do not arise in the algebraic group case.

It is easy to see as in the proof  \cite[Lemmas 2.1, 2.4]{LPS} that:

\begin{lemma} \label{tool1} Let $G$ be a classical group over $k$ with $n$
being the dimension of the natural module $V$.  The number of conjugacy classes of maximal subgroups
of $G$ in $ \cup_{H \in \mathcal{C}_i}$ is at most  $n +  3 \log n + 3$.
\end{lemma}

 We need to deal with the maximal subgroups in $\mathcal{S}_1$ and $\mathcal{S}_2$.
 Unfortunately in positive characteristic the results are rather weak and we do not have
 strong enough information about the subgroups in $\mathcal{S}_i, i =1,2$ to obtain
 good upper bounds on the number of conjugacy classes in $\mathcal{S}_i$.
  So now we assume that we are in characteristic $0$.
  
  By Theorem \ref{repeats}, it follows that the number of conjugacy classes of maximal
  subgroups in $\mathcal{S}_2(G)$ is $O(n^{\epsilon})$.
   
 The result for $\mathcal{S}_1$ we need has recently been obtained by  H\"as\"a \cite{hasa}.
 This is the only result in the paper that requires the classification of finite simple groups. 
 
 \begin{lemma}   The number of conjugacy classes of maximal subgroups in $\mathcal{S}_1(G)$
 is at most $O(n)$.
 \end{lemma}

 Now Theorem \ref{maximal}   for algebraic groups  follows -- for if the rank of $G$ is fixed,
 there is nothing to prove.  So we may assume that $G$ is a classical group in sufficiently
 large dimension and the remarks above immediately give the desired bound.
 
 The result for compact groups follows easily from the result on algebraic groups.

 \section{Dimensions of Fixed Spaces}
 
 In this section, we prove  Theorem \ref{gur-malle}.    Note that since $w(g)$ is constant
 on an open subvariety of $G$, it suffices to exhibit just one $g$ satisfying the conclusion
 (for a given $V$,  an open subvariety of a maximal torus has precisely the same weight
 spaces as the torus). 
 
 So let $G$ be a simple compact Lie group and $V$ a nontrivial irreducible finite dimensional
 module for $G$.   We require two lemmas.  The first is due  to Seitz \cite{Seitz} (and the result holds
 in all characteristics and for finite Chevalley groups).
 
 \begin{lemma} 
 \label{gary}
  Let $T$ be a maximal torus of $G$.  Then the dimension of any
 weight space for $T$ is at most $1 + (\dim V)/h$, where $h$ is the Coxeter number
 of $G$.
 \end{lemma}
 
 This proves the theorem if the rank is increasing.  So it suffices to prove the theorem for
 a fixed $G$. 
    Daniel Goldstein \cite{gold}  proved Theorem \ref{gur-malle}
 in this case. 
   We give a different proof (Goldstein's proof gives somewhat  better bounds). 
 
 \begin{lemma} 
 \label{danny}
Let $G$ be a simple compact Lie group.  Let $V$ be an irreducible finite dimensional module for $G$.
Let $g \in G$ be a regular element of prime order $p$.  Then
$$
\limsup_{ \dim V \rightarrow \infty}   \frac{w_V(g)}{\dim V} = \frac{1}{p}.
$$
\end{lemma}

\begin{proof}   Let $\chi$ denote the character of $V$. 
 By the Weyl character formula,
$$|\chi(g^k)| \le \frac{|W|}{|1-e^{2\pi i/p}|^{|\Phi^+|}}$$
for $1\le k\le p-1$ where $\Phi^+$ is the set of positive roots of
the root system of $G$ and $W$ is the Weyl group.    
Let $\lambda$ be any irreducible character for $L: =\langle g \rangle$.
If $m(\lambda, \chi)$ is the multiplicity of $\lambda$ in the restriction
of $\chi$ to $L$, then
$$ p m(\lambda, \chi)
  =  \sum_{k=0}^{p-1} \chi(g^k)\lambda(g^{-k})
  \le \chi(1) + (p-1) \frac{|W|}{|1-e^{2\pi i/p}|^{|\Phi^+|}}.$$
  
 Note that $w_V(g) = \max_{\lambda} \{ m(\lambda, \chi)\} \ge  \chi(1)/p$.   Thus
 the result follows. 
\end{proof}

Since for any sufficiently large prime $p$ there are regular elements of order $p$ in $G$,
we see that there is some element  $g  \in G$ with $w(g) < \epsilon \dim V$ as long
as $\dim V$ is sufficiently large.  This completes the proof.

\end{document}